\documentclass{conm-p-l}
\usepackage{amssymb}

\usepackage{amsfonts}
\usepackage{eucal}
\usepackage{upgreek}
\begin{document}

\font\eightrm=cmr8
\font\eightit=cmti8
\font\eighttt=cmtt8
\font\tensans=cmss10
\def\emp{\text{\tensans\O}}
\def\tci
{\hbox{\hskip1.8pt$\rightarrow$\hskip-11.5pt$^{^{C^\infty}}$\hskip-1.3pt}}
\def\nft
{\hbox{$n$\hskip3pt$\equiv$\hskip4pt$5$\hskip4.4pt$($mod\hskip2pt$3)$}}
\def\bbR{\mathrm{I\!R}}
\def\rto{\bbR\hskip-.5pt^2}
\def\rtr{\bbR\hskip-.7pt^3}
\def\rfo{\bbR\hskip-.7pt^4}
\def\rn{\bbR^{\hskip-.6ptn}}
\def\mr{\bbR^{\hskip-.6ptm}}
\def\bbZ{\mathsf{Z\hskip-4ptZ}}
\def\bbRP{\text{\bf R}\text{\rm P}}
\def\bbC{{\mathchoice {\setbox0=\hbox{$\displaystyle\rm C$}\hbox{\hbox
to0pt{\kern0.4\wd0\vrule height0.9\ht0\hss}\box0}}
{\setbox0=\hbox{$\textstyle\rm C$}\hbox{\hbox
to0pt{\kern0.4\wd0\vrule height0.9\ht0\hss}\box0}}
{\setbox0=\hbox{$\scriptstyle\rm C$}\hbox{\hbox
to0pt{\kern0.4\wd0\vrule height0.9\ht0\hss}\box0}}
{\setbox0=\hbox{$\scriptscriptstyle\rm C$}\hbox{\hbox
to0pt{\kern0.4\wd0\vrule height0.9\ht0\hss}\box0}}}}
\def\bbQ{{\mathchoice {\setbox0=\hbox{$\displaystyle\rm Q$}\hbox{\raise
0.15\ht0\hbox to0pt{\kern0.4\wd0\vrule height0.8\ht0\hss}\box0}}
{\setbox0=\hbox{$\textstyle\rm Q$}\hbox{\raise
0.15\ht0\hbox to0pt{\kern0.4\wd0\vrule height0.8\ht0\hss}\box0}}
{\setbox0=\hbox{$\scriptstyle\rm Q$}\hbox{\raise
0.15\ht0\hbox to0pt{\kern0.4\wd0\vrule height0.7\ht0\hss}\box0}}
{\setbox0=\hbox{$\scriptscriptstyle\rm Q$}\hbox{\raise
0.15\ht0\hbox to0pt{\kern0.4\wd0\vrule height0.7\ht0\hss}\box0}}}}
\def\bbQ{{\mathchoice {\setbox0=\hbox{$\displaystyle\rm Q$}\hbox{\raise
0.15\ht0\hbox to0pt{\kern0.4\wd0\vrule height0.8\ht0\hss}\box0}}
{\setbox0=\hbox{$\textstyle\rm Q$}\hbox{\raise
0.15\ht0\hbox to0pt{\kern0.4\wd0\vrule height0.8\ht0\hss}\box0}}
{\setbox0=\hbox{$\scriptstyle\rm Q$}\hbox{\raise
0.15\ht0\hbox to0pt{\kern0.4\wd0\vrule height0.7\ht0\hss}\box0}}
{\setbox0=\hbox{$\scriptscriptstyle\rm Q$}\hbox{\raise
0.15\ht0\hbox to0pt{\kern0.4\wd0\vrule height0.7\ht0\hss}\box0}}}}
\def\dimr{\dim_{\hskip.4pt\bbR\hskip-1.2pt}\w}
\def\dimc{\dim_{\hskip.4pt\bbC\hskip-1.2pt}\w}
\def\aff{\mathrm{A\hn f\hh f}\hs}
\def\cx{C\hskip-2pt_x\w}
\def\cy{C\hskip-2pt_y\w}
\def\cz{C\hskip-2pt_z\w}
\def\hyp{\hskip.5pt\vbox
{\hbox{\vrule width3ptheight0.5ptdepth0pt}\vskip2.2pt}\hskip.5pt}
\def\er{r}
\def\es{s}
\def\df{d\hskip-.8ptf}
\def\dz{\mathcal{D}}
\def\dzp{\dz^\perp}
\def\fv{\mathcal{F}}
\def\gr{\mathcal{G}}
\def\fvp{\fv_{\nrmh p}}
\def\wv{\mathcal{W}}
\def\vt{\mathcal{P}}
\def\tv{\mathcal{T}}
\def\vr{\mathcal{V}}
\def\xs{J}
\def\xl{S}
\def\cs{\mathcal{B}}
\def\zy{\mathcal{Z}}
\def\vtx{\vt_{\nh x}}
\def\fh{f}
\def\g{\mathtt{g}}
\def\rc{\theta}
\def\jm{\mathcal{I}}
\def\ke{\mathcal{K}}
\def\xc{\mathcal{X}_c}
\def\lz{\mathcal{L}}
\def\dla{\mathcal{D}_{\hskip-2ptL}^*}
\def\Lie{\pounds}
\def\lv{\Lie\hskip-1.2pt_v\w}
\def\lo{\lz_0}
\def\xe{\mathcal{E}}
\def\eo{\xe_0}
\def\lsq{\mathsf{[}}
\def\rsq{\mathsf{]}}
\def\hga{\hskip2.3pt\widehat{\hskip-2.3pt\gamma\hskip-2pt}\hskip2pt}
\def\hm{\hskip1.9pt\widehat{\hskip-1.9ptM\hskip-.2pt}\hskip.2pt}
\def\hg{\hskip.9pt\widehat{\hskip-.9pt\g\hskip-.9pt}\hskip.9pt}
\def\hna{\hskip.2pt\widehat{\hskip-.2pt\nabla\hskip-1.6pt}\hskip1.6pt}
\def\hdz{\hskip.9pt\widehat{\hskip-.9pt\dz\hskip-.9pt}\hskip.9pt}
\def\hdp{\hskip.9pt\widehat{\hskip-.9pt\dz\hskip-.9pt}\hskip.9pt^\perp}
\def\hmt{\hskip1.9pt\widehat{\hskip-1.9ptM\hskip-.5pt}_t}
\def\hmz{\hskip1.9pt\widehat{\hskip-1.9ptM\hskip-.5pt}_0}
\def\hmp{\hskip1.9pt\widehat{\hskip-1.9ptM\hskip-.5pt}_p}
\def\hk{\hskip1.5pt\widehat{\hskip-1.5ptK\hskip-.5pt}\hskip.5pt}
\def\hq{\hskip1.5pt\widehat{\hskip-1.5ptQ\hskip-.5pt}\hskip.5pt}
\def\txm{{T\hskip-2.9pt_x\w\hn M}}
\def\tyhm{{T\hskip-3.5pt_y\w\hm}}
\def\q{q}
\def\bq{\hat q}
\def\p{p}
\def\w{^{\phantom i}}
\def\x{u}
\def\y{v}
\def\vp{{\tau\hskip-4.55pt\iota\hskip.6pt}}
\def\evp{{\tau\hskip-3.55pt\iota\hskip.6pt}}
\def\vd{\vt\hh'}
\def\vdx{\vd{}\hskip-4.5pt_x}
\def\bz{b\hh}
\def\fe{F}
\def\fy{\phi}
\def\vl{\Lambda}
\def\hy{\mathcal{V}}
\def\vh{h}
\def\bc{C}
\def\mv{V}
\def\vo{V_{\nnh0}}
\def\ao{A_0}
\def\bo{B_0}
\def\uv{\mathcal{U}}
\def\sv{\mathcal{S}}
\def\svp{\sv_p}
\def\xv{\mathcal{X}}
\def\xvp{\xv_p}
\def\yv{\mathcal{Y}}
\def\yvp{\yv_p}
\def\zv{\mathcal{Z}}
\def\zvp{\zv_p}
\def\cv{\mathcal{C}}
\def\dy{\mathcal{D}}
\def\nv{\mathcal{N}}
\def\iv{\mathcal{I}}
\def\gkp{\Sigma}
\def\ret{\sigma}
\def\taw{\uptau}
\def\hs{\hskip.7pt}
\def\hh{\hskip.4pt}
\def\hn{\hskip-.4pt}
\def\nh{\hskip-.7pt}
\def\nnh{\hskip-1pt}
\def\hrz{^{\hskip.5pt\text{\rm hrz}}}
\def\vrt{^{\hskip.2pt\text{\rm vrt}}}
\def\vt{\varTheta}
\def\mtr{\Theta}
\def\op{\varTheta}
\def\vg{\varGamma}
\def\my{\mu}
\def\ny{\nu}
\def\gy{\lambda}
\def\lp{\lambda}
\def\ax{\alpha}
\def\lf{\widetilde{\lp}}
\def\bx{\beta}
\def\ay{a}
\def\by{b}
\def\gp{\mathrm{G}}
\def\hp{\mathrm{H}}
\def\kp{\mathrm{K}}
\def\gm{\gamma}
\def\Gm{\Gamma}
\def\Lm{\Lambda}
\def\Dt{\Delta}
\def\dg{\Delta}
\def\sj{\sigma}
\def\lg{\langle}
\def\rg{\rangle}
\def\lr{\langle\hh\cdot\hs,\hn\cdot\hh\rangle}
\def\vs{vector space}
\def\rvs{real vector space}
\def\vf{vector field}
\def\tf{tensor field}
\def\tvn{the vertical distribution}
\def\dn{distribution}
\def\pt{point}
\def\tc{tor\-sion\-free connection}
\def\ea{equi\-af\-fine}
\def\rt{Ric\-ci tensor}
\def\pde{partial differential equation}
\def\pf{projectively flat}
\def\pfs{projectively flat surface}
\def\pfc{projectively flat connection}
\def\pftc{projectively flat tor\-sion\-free connection}
\def\su{surface}
\def\sco{simply connected}
\def\psr{pseu\-\hbox{do\hs-}Riem\-ann\-i\-an}
\def\inv{-in\-var\-i\-ant}
\def\trinv{trans\-la\-tion\inv}
\def\feo{dif\-feo\-mor\-phism}
\def\feic{dif\-feo\-mor\-phic}
\def\feicly{dif\-feo\-mor\-phi\-cal\-ly}
\def\Feicly{Dif-feo\-mor\-phi\-cal\-ly}
\def\diml{-di\-men\-sion\-al}
\def\prl{-par\-al\-lel}
\def\skc{skew-sym\-met\-ric}
\def\sky{skew-sym\-me\-try}
\def\Sky{Skew-sym\-me\-try}
\def\dbly{-dif\-fer\-en\-ti\-a\-bly}
\def\cf{con\-for\-mal\-ly flat}
\def\ls{locally symmetric}
\def\ecs{essentially con\-for\-mal\-ly symmetric}
\def\rr{Ric\-ci-re\-cur\-rent}
\def\kf{Kil\=ling field}
\def\om{\omega}
\def\vol{\varOmega}
\def\og{\varOmega\hs}
\def\dv{\delta}
\def\ve{\varepsilon}
\def\zt{\zeta}
\def\kx{\kappa}
\def\mf{manifold}
\def\mfd{-man\-i\-fold}
\def\bmf{base manifold}
\def\bd{bundle}
\def\tbd{tangent bundle}
\def\ctb{cotangent bundle}
\def\bp{bundle projection}
\def\prc{pseu\-\hbox{do\hs-}Riem\-ann\-i\-an metric}
\def\prd{pseu\-\hbox{do\hs-}Riem\-ann\-i\-an manifold}
\def\Prd{pseu\-\hbox{do\hs-}Riem\-ann\-i\-an manifold}
\def\npd{null parallel distribution}
\def\pj{-pro\-ject\-a\-ble}
\def\pd{-pro\-ject\-ed}
\def\lcc{Le\-vi-Ci\-vi\-ta connection}
\def\vb{vector bundle}
\def\vbm{vec\-tor-bun\-dle morphism}
\def\kerd{\text{\rm Ker}\hskip2.7ptd}
\def\ro{\rho}
\def\sy{\sigma}
\def\ts{total space}
\def\pmb{\pi}
\def\pl{\partial}
\def\cro{\overline{\hskip-2pt\pl}}
\def\ddb{\pl\hskip1.7pt\cro}
\def\dbd{\cro\hskip-.3pt\pl}
\newtheorem{theorem}{Theorem}[section] 
\newtheorem{proposition}[theorem]{Proposition} 
\newtheorem{lemma}[theorem]{Lemma} 
\newtheorem{corollary}[theorem]{Corollary} 
  
\theoremstyle{definition} 
  
\newtheorem{defn}[theorem]{Definition} 
\newtheorem{notation}[theorem]{Notation} 
\newtheorem{example}[theorem]{Example} 
\newtheorem{conj}[theorem]{Conjecture} 
\newtheorem{prob}[theorem]{Problem} 
  
\theoremstyle{remark} 
  
\newtheorem{remark}[theorem]{Remark}

\title[Corrections of misstatements in several papers]{Corrections of minor
misstatements in several papers\\ on ECS manifolds}
\author[A. Derdzinski]{Andrzej Derdzinski} 
\address{Department of Mathematics, The Ohio State University, 
Columbus, OH 43210, USA} 
\email{andrzej@math.ohio-state.edu} 
\author[I.\ Terek]{Ivo Terek} 
\address{Department of Mathematics, The Ohio State University, 
Columbus, OH 43210, USA} 
\email{terekcouto.1@osu.edu} 
\def\leftmark{A.\ Derdzinski \&\ I.\ Terek}
\def\rightmark{Corrections of misstatements in several papers}

\begin{abstract}
In the paper \cite{derdzinski-roter-09} the value of the invariant denoted
by $\,d\,$ and often called "\nnh rank" was misstated for a narrow class of
examples of ECS manifolds: they were identified as having $\,d=1\,$ instead
of the correct value $\,d=2$. The error was repeated, by citing 
\cite{derdzinski-roter-09}, in \cite{derdzinski}, \cite{derdzinski-roter-07}
and \cite{derdzinski-terek-ne} -- \cite{derdzinski-terek-tc}. We briefly
describe the class in question, explain why it has $\,d=2$, and list the
required corrections of the affected papers.
\end{abstract}

\maketitle

\setcounter{section}{0}
\setcounter{theorem}{0}
\renewcommand{\thetheorem}{\Alph{theorem}}
\setcounter{equation}{0}
\section{Why, in some cases, $\,d=2$}
Given a pseu\-do\hs-Riem\-ann\-i\-an manifold $\,(M\nh,g)\,$ with parallel
Weyl tensor $\,W\nnh$, we denote by $\,d\,$ the rank of the
parallel sub\-bun\-dle
of $\,T^*\hskip-2ptM\nh$, the sections of which are the $\,1$-forms
$\,\xi\,$ such that $\,W(v,v'\nh,\,\cdot\,,\,\cdot\,)\wedge\xi=0\,$ for all 
vector fields $\,v,v'$ or, in other words, the $\,2$-form
$\,W(v,v'\nh,\,\cdot\,,\,\cdot\,)\,$ is $\,\wedge$-di\-vis\-i\-ble by
$\,\xi\,$ at points where $\,\xi\ne0$.

In a class of such manifolds, of dimensions $\,n\ge4$,  
constructed by Wi\-told Roter in 1974, $\,M\,$ is an open set in 
$\,\bbR\nh^n$ and, with 
$\,\lambda,\mu,\nu\,$ ranging over 
$\,2,\dots,n-1$, the pos\-si\-bly-non\-ze\-ro components of the metric
$\,g\,$ are those algebraically related to 
$\,g_{1\nh1}\w=[f(x^1)g_{\lambda\mu}\w+a_{\lambda\mu}\w]x^\lambda x^\mu\nh$,
$\,g_{1n}\w=g_{n1}\w=1$, and $\,g_{\lambda\mu}\w$. (We call them the {\it
essential components\/} of $\,g$.) Here $\,f\,$ is a function of $\,x^1$
and $\,[g_{\lambda\mu}\w],\,[a_{\lambda\mu}\w]\,$ are 
$\,(n-2)\times(n-2)$ symmetric matrices of constants such that
$\,\det\hs[g_{\lambda\mu}\w]\ne0\,$ and 
$\,g^{\lambda\mu}a_{\lambda\mu}\w=0$, while
$\,[a_{\lambda\mu}\w]\ne0$, where, as usual, 
$\,[g^{\lambda\mu}]=[g_{\lambda\mu}\w]^{-\nnh1}\nh$.

The reciprocal metric, Le\-vi-Ci\-vi\-ta connection, and the curvature,
Ric\-ci and Weyl tensors then have the essential components
$\,g^{1n}\nh=1$, $\,g^{nn}\nh=-g_{1\nh1}$, $\,g^{\lambda\mu}\nh$,
$\,\vg_{\hskip-2.7pt1\nh1}^{\hs\lambda}
=-g^{\lambda\mu}\partial\nh_\mu\w g_{1\nh1}\w/2$, 
$\,\vg_{\hskip-2.7pt1\nh1}^{\hs n}=\partial\nh_1\w g_{1\nh1}\w/2$,
$\,\vg_{\hskip-2.7pt1\nh\lambda}^{\hs n}=\partial\nh_\lambda\w g_{1\nh1}\w/2$,
$\,R_{1\lambda\mu1}\w=f(x^1)g_{\lambda\mu}\w+a_{\lambda\mu}\w$, 
$\,R_{1\nh1}\w=(2-n)f(x^1)$, and
$\,W\hskip-2.3pt_{1\lambda\mu1}\w=a_{\lambda\mu}\w$, 
the sign convention being 
$\,R\hn_{ij}\w=R\hn_{ikj}\w{}^k\nh$. All $\,2$-forms
$\,\zeta=W(v,v'\nh,\,\cdot\,,\,\cdot\,)$, thus lie, at each point,  
in the span of $\,\zeta_\lambda\w
=W(\partial\nh_1\w,\partial\nh_\lambda\w,\,\cdot\,,\,\cdot\,)
=a_{\lambda\mu}\w dx^\mu\wedge dx^1\nh$, where $\,\lambda=2,\dots,n-1$. As 
$\,[a_{\lambda\mu}\w]\ne0$, the condition 
$\,0=\zeta_\lambda\w\wedge\xi
=-a_{\lambda\mu}\w dx^1\wedge dx^\mu\wedge(\xi_\nu\w dx^\nu\nh
+\xi_n\w dx^n)$, for all $\,\lambda$, is equivalent to $\,\xi_n\w=0$ and 
$\,a_{\lambda\mu}\w\xi_\nu\w=a_{\lambda\nu}\w\xi_\mu\w$. Depending on
whether $\,\mathrm{rank}\hs[a_{\lambda\mu}\w]\ge2\,$ or 
$\,\mathrm{rank}\hs[a_{\lambda\mu}\w]=1$, this amounts to
$\,\xi_2\w=\ldots=\xi_n\w=0\,$ or, respectively, to requiring that
$\,\xi_n\w=0$ and 
$\,(\xi_2\w,\dots,\xi_{n-1}\w)$ be a multiple of every nonzero row of
$\,[a_{\lambda\mu}\w]$. Consequently -- with $\,\xi_1\w$ completely arbitrary
-- {\it in the former case\/ $\,d=1\,$ and in the latter\/ $\,d=2$}.

As $\,g_{in}\w=\delta_{1i}\w$ and $\,\vg_{\hskip-2.7ptin}^{\hs j}=0\,$ for all
$\,i,j=1,\dots,n$, the coordinate vector field $\,\partial\nh_n\w$ is null and
parallel, while $\,g(\partial\nh_n\w,\,\cdot\,)=dx^1\nh$. Thus, with $\,dx^1$
obviously being a section of the parallel rank $\,d\,$ sub\-bun\-dle of
$\,T^*\hskip-2ptM\,$ mentioned above, the rank $\,d\,$ sub\-bun\-dle of 
$\,T\nh M\,$ corresponding to it under the $\,g$-induced iso\-mor\-phism
$\,T^*\hskip-2ptM\to T\nh M$ contains the 
{\it one-di\-men\-sion\-al null parallel distribution\/} $\,\mathcal{P}$
{\it spanned by\/} $\,\partial\nh_n\w$.

\section{List of corrections: the two older papers}
The corrections needed in \cite{derdzinski-roter-09} and 
\cite{derdzinski-roter-07} are nearly identical.

In \cite[Theorem 4.1]{derdzinski-roter-09}, {\it `\hh has\/
$\,d=1$\hskip-2.4pt'} should read {\it `\hh has\/ $\,d=1\,$ or\/ $\,d=2\,$
depending on whether\/ $\,\mathrm{rank}\hs A\ge2\,$ or\/
$\,\mathrm{rank}\hs A=1$\hskip-2.4pt'}, and 
{\it `some such data\/  
$\,I\hskip-.7pt,f,n,V\hskip-.7pt,\langle\,,\rangle,A$\hskip-1.7pt'} needs to be
replaced by {\it `some such data\/  
$\,I\hskip-.7pt,f,n,V\hskip-.7pt,\langle\,,\rangle,A\,$ with\/
$\,\mathrm{rank}\hs A\ge2$\hskip-1.7pt'}.

In \cite[Theorem 3.1]{derdzinski-roter-07}: instead of 
{\it `\hh has parallel Weyl tensor and\/ $\,d=1$\hskip-2.4pt'} read\break 
{\it `\hh has parallel Weyl tensor and\/ $\,d=1\,$ or\/ $\,d=2\,$ depending on
whether\/ $\,\mathrm{rank}\hs A\ge2$ or\/
$\,\mathrm{rank}\hs A=1$\hskip-2.4pt'}, 
while {\it `some such\/  
$\,I\hskip-.7pt,f,n,V\hskip-.7pt,\langle\,,\rangle,A$\hskip-1.7pt'} should
similarly be replaced by {\it `some such\/  
$\,I\hskip-.7pt,f,n,V\hskip-.7pt,\langle\,,\rangle,A\,$ with\/
$\,\mathrm{rank}\hs A\ge2$\hskip-2.3pt'}.

\section{List of corrections: \cite{derdzinski-terek-ne}}
The corrected version of lines 5 -- 3 before Theorem A is

\noindent\hbox{`locally-homogeneous\hn\ compact\hn\ rank-two\hn\ ECS\hn\
manifolds\hn\ of\hn\ all\hn\ 
{\it odd\/}\hn\ dimensions\hn\ $n\nh\ge\nh5$.}

\noindent It still remains an open question whether a compact ECS
  manifold
may be of dimension four.'

Line 7 before formula (2.2) should read 
`The metric (2.1) turns $\,\widehat{M}\,$ into an ECS manifold
[5, Theorem 4.1], of rank one or two depending on
whether\/ $\,\mathrm{rank}\hs A\ge2$ or\/
$\,\mathrm{rank}\hs A=1$.'

In the third line of Remark 2.1, 
`are as before' needs to be replaced by `are as before, with
$\,\mathrm{rank}\hs A\ge2$'.

In the line preceding formula (3.1), 
`\hs with the metric (2.1)' should read `\hs with the metric (2.1), where
$\,\mathrm{rank}\hs A\ge2$'.

\section{List of corrections: \cite{derdzinski-terek-ro}}
Abstract: lines 4 -- 5 should read `Ol\-szak, and a Ric\-ci-re\-cur\-rent ECS
manifold may be called translational or di\-la\-tion\-al, depending on whether
the holonomy group of a natural flat connection in a specific line sub\-bundle
of the Olszak distribution is finite'.

Abstract: lines 7 -- 9 need to be replaced with `structure of the Weyl
tensor. Various examples of compact ECS manifolds are known: rank-one
translational ones (both generic and nongeneric) in every dimension $\,n\ge5$,
as well as odd-di\-men\-sion\-al rank-two nongeneric dilational ones, some of
which are locally homogeneous. As we show,'.

P.\ 70, lines 5 -- 8 should read `is locally homogeneous. Quite recently [8]
we constructed dilational type compact rank-two ECS manifolds, including
locally-homogeneous ones, in all odd dimensions $\,n\ge5$. It remains an open
question whether a compact ECS manifold may be of dimension four'.

The line before formula (4.3) has to be replaced by
`It is well known [4, Theorem 4.1] that (4.2) is an ECS manifold, of
rank one if $\,\mathrm{rank}\hs A>1$, and of rank two when
$\,\mathrm{rank}\hs A=1$. To'.

\section{List of corrections: \cite{derdzinski-terek-ms}}
Lines 7 -- 10 before Theorem A should read 'have recently found [26] rank-two
compact dilational examples in all odd dimensions $\,n\ge5$, including locally
homogeneous ones. They are all nongeneric and incomplete. 
It is still not known whether a compact ECS manifold can be four-dimensional.'

The second line after formula (6.2) has to be replaced by 
`By [4, Theorem 4.1] $\,(\widehat M\nh,\widehat g)\,$ is an ECS manifold, of 
rank one if $\,\mathrm{rank}\hs A>1$, and of rank two when 
$\,\mathrm{rank}\hs A=1$. Calling the manifolds (6.2)'

\section{List of corrections: \cite{derdzinski-terek-tc}}
Abstract, line 4: `are of rank one' should read `are of rank one (except for
some odd-di\-men\-sion\-al rank-two examples, found recently)'

Three lines before Theorem A: `compact locally homogeneous ECS manifolds'
needs to be replaced with `compact locally homogeneous rank-two ECS manifolds'.

Lines 3 -- 4 before Theorem B should read `We do not know whether the final
clause of Theorem A is non-vacuous.'

In formula (5.1), replace `and a non-zero, trace\-less,
$\,\langle\cdot,\cdot\rangle$-self-ad\-joint linear en\-do\-mor\-phism $\,A\,$ of $\,V$\hh' by `and a trace\-less
$\,\langle\cdot,\cdot\rangle$-self-ad\-joint linear en\-do\-mor\-phism $\,A\,$ of $\,V\hs$ with $\,\mathrm{rank}\hs A\ge2$'.

\section{List of corrections: \cite{derdzinski}}\label{cw}
Abstract: lines 7 -- 9 should read `compact ECS manifolds [4,\ 6],
representing every dimension $\,n\ge5$, are of rank $\,1$, except for
some recently found [7] odd-di\-men\-sion\-al locally homogeneous examples of 
rank $\,2$'.

Section 2: lines 3 -- 4 need to be replaced with `homogeneous examples of
rank two (which are necessarily incomplete) were exhibited in all odd
dimensions $\,n\ge5$. It is an open question whether a compact ECS manifold
may have dimension $\,4$'.

The second line before Theorem 3.2 should read 
`locally homogeneous compact rank-two ECS manifolds, constructed in [7], are
bundles over the circle'.

The corrected version of the line following formula (7.1) is
`defines a rank-one ECS metric if $\,f\,$ is nonconstant, 
$\,\det\hh[g_{ij}\w]\ne0=g\hh^{ij}a_{ij}\w$ and
$\,\mathrm{rank}\hs[a_{ij}\w]\ge2$'.

\end{document}